\documentclass[twoside]{article}
\usepackage[a4paper]{geometry}
\usepackage[latin1]{inputenc} 
\usepackage[T1]{fontenc} 
\usepackage{RR}
\usepackage{hyperref}
\usepackage{amsmath, amssymb}
\usepackage{algorithmic}
\usepackage{algorithm}

\RRNo{8265}
\RRdate{February 2013}
\RRauthor{
Olivier Coulaud\thanks{Inria Bordeaux-Sud Ouest, France}
\and
Luc Giraud\footnotemark[1]
\and
Pierre Ramet\thanks{Universit\'e de Bordeaux 1, France}
\and
Xavier Vasseur\thanks{CERFACS, France}
}
\authorhead{Coulaud \& Giraud \& Ramet \& Vasseur}
\RRtitle{Techniques de d\'eflation et d'augmentation dans les solveurs lin\'eaires de Krylov}
\RRetitle{Deflation and augmentation techniques in Krylov subspace methods for the solution
 of linear systems}
\titlehead{Deflation and augmentation techniques in Krylov lienar solvers}
\RRresume{
Dans ce rapport nous pr\'esentons des techniques de d\'eflation et d'augmentation qui ont \'et\'e d\'evelopp\'ees pour acc\'el\'erer
la convergence des m\'ethodes de Krylov pour la solution de syst\`emes d'\'equations lin\'eaires.
Nous passons en revue des approches pour des syst\`emes lin\'eaires dont les matrices sont non-hermitiennes, principalement
dans le contexte de la m\'ethode d'Arnoldi, et pour des matrices hermitiennes d\'efinies positives avec la m\'ethode du gradient conjugu\'e.
}
\RRabstract{
In this paper we present deflation and augmentation techniques that have been designed to accelerate the convergence of Krylov subspace methods
for the solution of linear systems of equations. 
We review numerical approaches both for linear systems with a non-Hermitian coefficient matrix, mainly within the Arnoldi framework, and for Hermitian
positive definite problems with the conjugate gradient method.
}
\RRmotcle{
Augmentation, D\'eflation, M\'ethodes de Krylov, Syst\`emes lin\'eaires d'\'equations, Pr\'econditionnement.
}
\RRkeyword{
Augmentation, Deflation, Krylov subspace methods, Linear systems of equations, Preconditioning.
 }
\RRprojets{HiePACS}
 \RCBordeaux 

\begin{document}
\makeRR   
\tableofcontents
\newpage

\newtheorem{definition}{Definition}
\newtheorem{lemma}{Lemma}

\newtheorem{Prop}{Proposition}
\newtheorem{Lem}{Lemma}
\newtheorem{Rq}{Remark}

\newcommand{\Complex}{\mathbb{C}}
\newcommand{\Prec}{{\cal{M}}}

\newcommand{\Hp}{P_{{W^{\perp_A}}}}
\newcommand{\argmin}{\mathop{\mathrm{argmin}}}
\newcommand{\norm}[1]{\Vert #1 \Vert}
\newcommand{\gcrodr}[1]{{#1^{}}}
\newcommand{\gmresdr}[1]{\widetilde{#1}}
\newcommand{\gmresdralg}[1]{#1}
\newcommand{\range}[1]{\ensuremath{\operatorname{Range}(#1)}}

\newcommand{\MDeflat}{M^{def}}
\newcommand{\MShift}{M^{coarse}}

\newcommand{\Qu}{Q_1}
\newcommand{\Qd}{Q_2}
\newcommand{\Qq}{Q_3}
\newcommand{\Pu}{P_1}
\newcommand{\Pd}{P_2}
\newcommand{\Pt}{P_3}
\newcommand{\Pq}{P_6}
\newcommand{\Pc}{P_4}
\newcommand{\Ps}{P_5}

\section{\label{sec:intr} Introduction}

The solution of linear systems  of the form $A x^\star = b$ plays a central role in many engineering and academic simulation codes.
Among the most widely used solution techniques are the iterative schemes based on Krylov subspace methods \cite{axel:94,erhe:11,saad:03,vors:03}.
Their main advantages are their ability to solve linear systems even if the matrix of the linear system is not
explicitly available and their capability to be ``easily'' parallelizable on large computing platforms.
In order to speed up the convergence of these solution techniques, Krylov subspace methods are almost always used in combination with 
preconditioning. That is, instead of solving directly $Ax^\star = b$, the linear system is transformed into an equivalent one, e.g.,  $M_1Ax^\star = M_1b$, referred to as left preconditioned system, that is expected to be more amenable to a  solution.
The definition of an efficient preconditioner $M_1$, that should be an good approximation of $A^{-1}$ in some sense,
is very much problem dependent and is consequently an extremely active research field. 
We can also consider other equivalent linear systems $A M_2 t^\star = b$ with $x^\star=M_2 t^\star$ (right preconditioner) or $M_1 A M_2 t^\star = M_1 b$ with $x^\star = M_2 t^\star$ (split preconditioner). We refer the reader to \cite{benz:02} for a detailed overview on preconditioning. 

There exist two complementary alternatives to speed up the convergence of the Krylov space, namely augmentation and deflation.
Roughly speaking, in augmentation techniques,  the search space in an enlarged Krylov space that is defined by
a direct sum of two subspaces. This search space $S_{\ell}$ (of dimension $\ell$) has the following form
\begin{equation}\label{eq:def S}
S_\ell = {\cal{K}}_m(A, b) \oplus {\cal{W}}
\end{equation}
where ${\cal{K}}_m(A, b)$ is a Krylov subspace of dimension $m$ generated by the matrix $A$ and the vector $b$ and ${\cal{W}}$ (of dimension $k$) is called the augmentation space. A typical goal of augmentation is to add information about the problem into the global search space $S_{\ell}$ that is only slowly revealed in the Krylov subspace itself. 

Alternatively, deflation is based on the use of a projection operator $P$ to decompose $x^\star$ as $x^\star = P x^\star + (I-P)x^\star$. The general idea is to select $P$ such that the solution of $PAx^\star = P b$, referred to as the deflated linear system, is more easily amenable to a solution by a Krylov subspace method than the original linear system $A x^\star = b$. The component $(I-P)x^\star$ can then be computed by solving a linear system of small dimension.

The purpose of this paper is to expose these two latter acceleration techniques that become increasingly popular. 
We refer the reader to \cite{gutk:12} for a recent excellent analysis of these methods together with detailed references and historical comments. Here augmentation and deflation are described in a framework where variable preconditioning can be used as it is nowadays customary when considering large scale linear systems \cite{nota:00,nova:08,sisz:07}.     
This paper is organized as follows. In Section~\ref{sec:kryl} we introduce some background on Krylov subspace methods with emphasis on the minimum residual norm approach for systems with a non-Hermitian coefficient matrix and the conjugate gradient method for the
solution of Hermitian positive definite problems. In Section~\ref{sec:nonh} we describe the augmentation and deflation techniques and their possible combination in the case of systems with non-Hermitian matrices with references to concrete applications. Similar exposure is performed in Section~\ref{sec:hpd} for Hermitian positive definite linear systems. Finally some concluding remarks and prospectives are drawn in Section~\ref{sec:conc}.

\section{\label{sec:kryl} Some background on Krylov subspace methods}

We briefly describe the basic properties of Krylov subspace methods for the solution of a linear 
system of equations of the form 
\begin{equation}
 A x^\star = b \label{eq:orig}
\end{equation}
where the nonsingular $n \times n$ coefficient matrix $A \in \mathbb{C}^{n\times n}$ is supposed to be either non-Hermitian or Hermitian 
positive definite and $b$ a given vector in $\mathbb{C}^{n}$. First we introduce the notation used throughout this paper.

\subsection{\label{sec:kryl:nota} Notation}

We denote the range of a matrix $A$ by $\range{A}$ and its nullspace as $Ker(A)$. 
We denote by $\| . \|$ the Euclidean norm, $I_k \in \mathbb{C}^{k \times k}$ the identity matrix of dimension $k$ and $0_{i \times j} \in \mathbb{C}^{i \times j}$ the zero rectangular matrix with $i$ rows and $j$ columns. ${}^T$ denotes the transpose operation, while ${}^H$ denotes the Hermitian transpose operation. Given a vector $d\in \Bbb{C}^{k}$ with components $d_i$, $D={\rm
diag}(d_1,\ldots, d_k)$ is the diagonal matrix $D\in \Bbb{C}^{k \times k}$ such that $D_{ii} = d_i$.  Given $Z_m = [z_1, \cdots, z_m] \in \Bbb{C}^{n \times m}$ we denote its $i$-th column as $z_i \in \mathbb{C}^{n},  (1 \le i \le m)$.
The vector $e_m \in \mathbb{R}^{m}$ denotes the $m$-th canonical basis vector of $\mathbb{R}^{m}$. 
We denote by $\kappa(A)$ the Euclidean condition number of $A$ that is defined by $\displaystyle \kappa(A)=\frac{\sigma_{\max}}{\sigma_{\min}}$ where
$\sigma_{\max}$ ($\sigma_{\min}$) is the largest (respectively smallest) singular value of $A$.
For Hermitian positive definite matrices, the condition number reduces to $\displaystyle \kappa(A)=\frac{\lambda_{\max}}{\lambda_{\min}}$ where
$\lambda_{\max}$ ($\lambda_{\min}$) is the largest (respectively smallest) eigenvalue of $A$.
Finally, throughout the paper for the sake of readability the integer subscript $\ell$ denotes the dimension of the search space.

\subsection{\label{sec:kryl:basi} Basic properties of Krylov subspace methods}
%
%
The Krylov subspace methods seek for the solution of Equation~\eqref{eq:orig} in a sequence of embedded spaces of increasing dimension ${\cal{K}}_{\ell}(A, b) = span (b, Ab, \cdots, A^{\ell-1}b)$. This is motivated~\cite{ipme:98} by the fact that for $\ell$ large enough these spaces contain the solution of the linear system~(\ref{eq:orig}).
If we denote by $m_A(x)$ the minimal polynomial associated with $A$, the Jordan decomposition of this polynomial writes
$m_A(t) = \prod_{i=1}^s (t - \lambda_i)^{m_i}$ where $(\lambda_1, \cdots, \lambda_s)$ are the distinct eigenvalues of $A$ and
$(m_1, \cdots, m_s)$ their indices in the Jordan form.
In a canonical form we also have $m_A(t) = \sum_{i=0}^m \alpha_i t^i$ with $m=\sum_{i=1}^s m_i$ and
$\alpha_0 = \prod_{i=1}^s (- \lambda_i)^{m_i} \neq 0$ since $A$ is nonsingular.
Consequently, $A^{-1} = - \alpha_0^{-1}\sum_{i=0}^{m-1} \alpha_{i+1} A^{i}$ that portrays $x^{\star}=A^{-1}b$ as a vector of the
Krylov space ${\cal{K}}_{m-1}(A, b)$.
This indicates that, in exact arithmetic, Krylov methods must converge in at most $m-1$ steps
or less if the right-hand side does not have components in all the eigendirections. 
This observation also gives some ideas on ways to speed-up the convergence of these methods.
As mentioned earlier preconditioning is a widely used approach that consists in transforming~\eqref{eq:orig} in an equivalent nonsingular system where the
preconditioned matrix has less~\cite{mugw:00} or better clustered eigenvalues (see~\cite{benz:02} and the references therein).

The rest of the paper is dedicated to an overview of proposed techniques for augmentation and deflation 
both for non-Hermitian and Hermitian positive definite problems.

\subsection{\label{sec:kryl:mini} Minimum residual Krylov subspace method}

In this section we focus on minimum residual norm subspace methods for the solution of linear systems with a non-Hermitian coefficient matrix. We refer the reader to \cite{saad:03,vors:03} for a general introduction to Krylov subspace methods and 
to \cite{sisz:07} for a recent overview on Krylov subspace methods; see also \cite{eier:01, eies:00} for an advanced analysis related to minimum residual norm Krylov subspace methods. \\

Augmented and deflated minimum residual norm Krylov subspace methods are usually characterized by a generalized Arnoldi relation introduced next.
\begin{definition}
\label{def:flex_arnoldi} Generalized Arnoldi relation. The minimum residual norm subspace methods investigated in this paper satisfies the following relation: 
\begin{equation}
  A Z_\ell = V_{\ell+1} \bar{H}_\ell \label{eq:def:flex_arnoldi}
\end{equation}
where $Z_\ell  \in \mathbb{C}^{n \times \ell}$, $V_{\ell+1}  \in \mathbb{C}^{n \times (\ell+1)}$ such that $V_{\ell+1}^H V_{\ell+1} = I_{\ell+1}$ and $\bar{H}_\ell \in \mathbb{C}^{(\ell+1) \times \ell}$. These methods compute
an approximation of the solution of (\ref{eq:orig}) in a $\ell$-dimensional affine space $x_0+Z_\ell \, y_{\ell}$ where $y_{\ell} \in \mathbb{C}^{\ell}$. In
certain cases, $\bar{H}_\ell $ is an upper Hessenberg matrix.
\end{definition}

We next introduce a minimum residual norm subspace method proposed by \linebreak[4] Saad~\cite{saad:93} since it is the basis for further developments related to augmented and deflated Krylov subspace methods of minimum residual norm type. 
This method named Flexible GMRES (FGMRES) was primarily introduced to allow variable preconditioning.
 We denote by $M_j$ the nonsingular matrix that represents the preconditioner at step $j$ of the method. Algorithm~\ref{alg:fgmres} depicts the FGMRES($\ell$) method where the dimension of the approximation subspace is not allowed to be larger than a prescribed dimension noted $\ell$. Starting from an initial guess $x_0 \in\mathbb{C}^{n}$, it is based on a generalized Arnoldi relation 
\begin{equation}
\label{eq:gmres-flex:arnoldi}
A Z_{\ell} = V_{\ell+1} \bar{H}_{\ell} \quad \mbox{with} \quad V_{\ell+1}^H V_{\ell+1} = I_{\ell+1},
\end{equation}
where $Z_{\ell} \in\mathbb{C}^{n \times \ell}$, $V_{\ell+1}\in\mathbb{C}^{n \times (\ell+1)}$ and the upper Hessenberg matrix $\bar{H}_{\ell}\in\mathbb{C}^{(\ell+1) \times \ell}$ are obtained from the Arnoldi procedure described in Algorithm \ref{alg:fgmres_comp}. An approximate solution  $x_{\ell} \in\mathbb{C}^n$ is then found by minimizing the residual norm $\norm{b-A(x_0+Z_{\ell} y)}$ over the space $x_0+ {\rm range}(Z_{\ell})$, the corresponding  residual being $r_{\ell}=b-Ax_{\ell}\in\mathbb{C}^n$ with $r_{\ell} \in {\rm range}(V_{{\ell}+1})$. With notation of Algorithm \ref{alg:fgmres} the current approximation $x_{\ell}$ can be written as 
\begin{equation}\label{eq:fgmresUpdate}
x_ {\ell}= x_0 + Z_{\ell} y^\star,
\end{equation}
whereas the residual $r_{\ell} = b -A x_{\ell}$ satisfies the Petrov-Galerkin orthogonality condition
\[
r_{\ell} \perp A~{\rm range}(Z_{\ell}).
\]
Hence, an optimality property similar to the one that defines GMRES is thus obtained~\cite{saad:03}. 
We note however that no general convergence results are available since the subspace of approximants $\range{Z_{\ell}}$ is no longer a standard Krylov subspace. 
We refer the reader to~\cite{saad:93,saad:03} for the analysis of the breakdown in FGMRES. 
Furthermore, as it can be seen in Equation~\eqref{eq:fgmresUpdate}, the update of the iterate $x_\ell$ requires to store the complete set of vectors $Z_\ell$ inducing a large memory footprint for large $\ell$.
In order to alleviate this memory requirement, a restarting strategy must be implemented as shown in Algorithm~\ref{alg:fgmres}.
The construction of a complete set of $Z_\ell$ is often name a cycle of the method and corresponds to one iteration of the loop in 
Algorithm~\ref{alg:fgmres}.

\begin{algorithm*}[!ht]
\caption{Flexible GMRES(${\ell}$)}
\label{alg:fgmres}
\begin{algorithmic}[1]
\STATE \textit{Initialization:}  
Choose ${\ell}>0$, $tol>0$, $x_0\in \mathbb{C}^n$. Let $r_0 = b - Ax_0$, $\beta = \norm{r_0} $, $c=[\beta, 0_{1\times {\ell}}]^T$ where $c\in \mathbb{C}^{{\ell}+1}$, $v_1=r_0/\beta$.

{\bf{Loop}}

\begin{minipage}[t]{0.05\linewidth}
\hspace{1cm} 
\end{minipage}
\hfill
\begin{minipage}[t]{0.95\linewidth}

\STATE \textit{Computation of $V_{{\ell}+1}$, $Z_{{\ell}}$ and $\bar{H}_{\ell}$ (see Algorithm \ref{alg:fgmres_comp}):}  
  Apply ${\ell}$ steps of the Arnoldi method with variable preconditioning ($z_j = M_j^{-1} v_j, 1 \le j \le {\ell}$) to obtain $V_{{\ell}+1} \in \mathbb{C}^{n \times ({\ell}+1)}$, 
$Z_\ell \in \mathbb{C}^{n \times {\ell}}$  and the upper Hessenberg matrix $\bar{H}_{\ell} \in \mathbb{C}^{({\ell}+1) \times {\ell}}$ such that
$$
      A Z_{\ell} = V_{{\ell}+1} \bar{H}_{\ell} \quad \mbox{with} \quad V_{{\ell}+1}^H  V_{{\ell}+1} = I_{{\ell}+1}.
$$

\STATE \textit{Minimum norm solution:} Compute the minimum norm solution $x_{\ell} \in \mathbb{C}^{n}$ in the affine space $x_0 + {\rm range}(Z_{\ell})$; that is, $x_{\ell} = x_0 + Z_ {\ell}y^\star$ where $\displaystyle y^\star = \argmin_{y \in \mathbb{C}^{{\ell}}} \norm{c - \bar{H}_{\ell} y}$.

\STATE \textit{Check the convergence criterion:} If $\norm{c - \bar{H}_{\ell} y^\star}/\norm{b} \leq tol$,  exit

\STATE \textit{Restarting:} Set $x_0= x_{\ell}$, $r_0 = b-Ax_0$, $\beta = \norm{r_0} $, $c=[\beta, 0_{1\times {\ell}}]^T$, $v_1=r_0/\beta$.

\end{minipage}

\vspace*{.1cm}
{\bf{End of loop}}

\end{algorithmic}
\end{algorithm*}

\begin{algorithm*}[!ht]
\caption{Arnoldi procedure: computation of $V_{{\ell}+1}$, $Z_{\ell}$ and $\bar{H}_{\ell}$}
\label{alg:fgmres_comp}
\begin{algorithmic}[1]
\FOR {$j=1, {\ell}$}
\STATE $z_j = M_j^{-1} v_j$
\STATE $s = A z_j$
\FOR {$i=1, j$}
\STATE     $h_{i,j} = v_i^H s$
\STATE     $s = s - h_{i,j} v_i$
\ENDFOR
\STATE $h_{i+1,j} = \norm{s}$, $v_{j+1} = s/h_{i+1,j}$
\ENDFOR
\STATE Define $Z_{\ell} = [z_1, \cdots, z_{\ell}]$, $V_{{\ell}+1} = [v_1, \cdots, v_{{\ell}+1}]$, $\bar{H}_ {\ell}= \{h_{i,j}\}_{1\le i \le {\ell}+1, 1\le j \le {\ell}}$
\end{algorithmic}
\end{algorithm*}

When the preconditioner is  constant, FGMRES(${\ell}$) reduces to right-preconditioned GMRES($\ell$) whose 
convergence properties are discussed in \cite[Chapter 6]{saad:03}.

\subsection{\label{sec:kryl:cg} Conjugate gradient method}
The conjugate gradient~\cite{hest:52} is the method of choice for Hermitian Positive Definite (HPD) linear systems.
In a shortcut, it relies on an Arnoldi like relation (namely a Lanczos relation~\cite{lanc:52}) similar to Equation~\eqref{eq:gmres-flex:arnoldi} with $Z_\ell = V_\ell$ (case of no preconditioning)
and a Ritz-Galerkin condition $r_\ell = b - A x_\ell \perp  {\cal{K}}_{\ell}(A, r_0)$.
At each iteration $x_\ell = x_0 + V_\ell y^{\star}$ is computed via the solution of the small linear system $H_\ell y^{\star} = \norm{r_0} (1, 0_{1\times (\ell-1)})^T$,
 where $H_\ell =  V_\ell^H A V_\ell$ is the square leading part of $\bar{H}_\ell$.
Since $A$ is Hermitian, $H_\ell$ is also Hermitian.
Furthermore, its structure is upper Hessenberg that, combined with the Hermitian property, implies that $H_\ell$ is tridiagonal HPD.
The first consequence of this structure of $H_\ell$ is that the orthogonalization of $V_\ell$ can be performed cheaply with a three term recurrence.
The second consequence is that a $LU$
factorization of $H_\ell$ can be incrementally computed and this factorization without pivoting is known to be stable for positive definite matrices.
The conjugate gradient method is a very elegant, sophisticated and powerful algorithm that exploits nicely all the above mentioned properties.
It can be implemented through short recurrences that do not require  to store the complete set of vector $V_\ell$ leading to a very low memory consumption.
Furthermore, the conjugate gradient enjoys a unique minimum norm property on the forward error that reads
$x_\ell = \argmin_{x \in x_0 + {\cal{K}}_{\ell}(A, r_0)} \norm{x - x^\star}_A$ where $x^\star$ denotes the exact solution and $\norm{ \cdot}_A$ is the norm associated with $A$.
 In addition, it exists an upper bound on its
convergence rate that reads ($\ell \ge 1$)
\begin{equation}\label{eq:CGCvg}
\displaystyle \norm{x_\ell - x^{\star}}_A \leq 
2\cdot \left ( \frac{ \sqrt{\kappa(A)} - 1}{\sqrt{\kappa(A)} + 1} \right)^\ell
    \norm{x_0 - x^{\star}}_A .
\end{equation}
We refer to~\cite{meur:06,saad:03,vors:03} for an exhaustive and detailed exposure of CG and to~\cite{gool:89} for a nice description of its history.

\section{\label{sec:nonh} Non-Hermitian matrices}

In this section we detail augmentation and deflation techniques in Krylov subspace methods when the coefficient matrix $A$ is non-Hermitian. We specifically focus on minimum residual norm subspace methods and assume that a generalized Arnoldi relation (\ref{eq:def:flex_arnoldi}) holds. We denote by $x_0$, $r_0=b-A x_0$ the initial guess and residual vector respectively, and by $V_{\ell+1}$, $\bar{H}_{\ell}$ and $Z_{\ell}$ the matrices involved in this relation. With notation of Algorithm 
\ref{alg:fgmres} $r_0$ can be expressed as $r_0 = V_{\ell+1} (c - \bar{H}_{\ell} y^{\star})$.

\subsection{\label{sec:nonh:aug} Augmented Krylov subspace methods}

  We next discuss two possibilities to select the augmentation space and analyze the corresponding Krylov subspace methods. 

\subsubsection{\label{sec:nonh:aug:arbr} Augmentation with an arbitrary subspace}

Given a basis $W=[w_1, \cdots, w_k]$ of an augmentation subspace ${\cal{W}}$, a slight modification in the Arnoldi procedure (Algorithm \ref{alg:fgmres_comp}) is used to deduce an 
orthogonal basis of $S_\ell$ defined in~\eqref{eq:def S} (see~\cite{chsa:97}). 
It consists of defining $z_j$ (line 2 of Algorithm~\ref{alg:fgmres_comp}) now as
\[
z_j = M_j^{-1} v_j \quad (1 \le j \le m) \; \mbox{\rm and }\;  z_j = M_j^{-1} w_{j-m} \quad (m < j \le m+k).
\]
With this definition we  finally obtain the generalized Arnoldi relation 
\[
A Z_{m+k} = V_{m+k+1} \bar{H}_{m+k}
\]
where
\begin{eqnarray}
Z_{m+k}   &=& [M_1^{-1} v_1, M_2^{-1} v_2, \cdots, M_{m+1}^{-1} w_1, M_{m+2}^{-1} w_2, \cdots, M_{m+k+1}^{-1} w_{k}], \\
V_{m+k+1} &=& [v_1, v_2, \cdots, v_{m+k+1}],
\end{eqnarray}
and $\bar{H}_{m+k}$ is a $(m+k+1) \times (m+k)$ upper Hessenberg matrix. Thus the residual minimization property 
is then deduced similarly as in FGMRES \cite{saad:93}. Hence, the approximate solution from the affine 
space $x_0 + {\rm range}(Z_{m+k})$ can be written as
\[
x_{m+k} = x_0 + Z_{m+k} y^{\star}
\]
with $y^{\star} \in \mathbb{C}^{(m+k)} $ solution of the residual norm minimization problem
\[
 y^\star = \argmin_{y \in \mathbb{C}^{(m+k)}} \norm{||r_0|| e_1 - \bar{H}_{m+k} y}
\]
(with $e_1$ designing here the first canonical  vector of $\mathbb{R}^{(m+k+1)}$). 
In case of constant right preconditioning the main important property is that if any vector $w_j$ is the solution of $AM^{-1}w_j = v_i, 1 \le i \le m$,  then in general the
exact solution of the original system~(\ref{eq:orig}) can be extracted from $S_\ell$; see, e.g., 
\cite[Proposition 2.1]{saad:97}. We refer the reader to \cite{chsa:97} for a discussion of possible choices for the augmented subspace ${\cal{W}}$. Vectors obtained with either different iterative methods or with different preconditioners can  be incorporated in $Z_{m}$ quite 
easily. In block Krylov subspace methods we also mention that ${\cal{W}}$ consists of the sum of a few other Krylov subspaces generated with the same matrix but with different right-hand sides;
 see~\cite{chsa:97} for a discussion and numerical experiments on academic problems. A popular idea is to choose ${\cal{W}}$ as an approximate invariant subspace associated with a specific part of the spectrum of $A$ or $AM^{-1}$ in case of fixed preconditioning. This is discussed next. 

\subsubsection{\label{sec:nonh:aug:app} Augmentation with approximate invariant subspace}

A typical goal of augmentation is to add information about the problem into the search space that is only slowly revealed in the Krylov subspace itself. It is often known that eigenvalues of the (preconditioned) operator close to zero tend to slow down the convergence rate of the Krylov subspace methods~\cite{chsa:97}. Hence, augmentation based on approximate invariant subspaces made of eigenvectors corresponding to small in modulus eigenvalues of the (preconditioned) operator has been proposed; see, e.g.,\cite{morg:95,morg:00,morg:02,saad:97} and references therein.

\paragraph{\label{sec:nonh:aug:app:hr} Harmonic Ritz information}

In \cite{morg:02} Morgan has suggested to select ${\cal{W}}$ as an approximate invariant subspace and to update this subspace at
the end of each cycle. Approximate spectral information is then required to define the augmentation space. This is usually obtained by computing harmonic Ritz pairs of $A$ with respect to a certain subspace \cite{chsa:97,morg:02}. We present here a definition of a harmonic Ritz pair as given in \cite{papv:95,slvdv:1996}. 

\begin{definition}
\label{def:hr} Harmonic Ritz pair. Consider a subspace $~{\cal{U}}$ of $ \mathbb{C}^{n}$. 
Given $B \in\mathbb{C}^{n \times n}$, $\theta \in \mathbb{C}$ and $y \in \cal {U}$, $(\theta, y)$ is a harmonic Ritz pair of $B$ with respect to $\cal {U}$ if and only if
\[
B y - \theta \, y \perp B \, {\cal {U}}
\]
or equivalently, for the canonical scalar product,
\[
\forall w \in \range{B \, {\cal {U}}} \quad w^H \, ( B y - \theta \, y) = 0.
\]
We call the vector $y$ a harmonic Ritz vector associated with the harmonic Ritz value $\theta$. 
\end{definition}

Based on the generalized Arnoldi relation (\ref{eq:def:flex_arnoldi}), the augmentation procedure proposed in \cite[Proposition 1]{ggpv:10} relies on the use of $k$ harmonic Ritz vectors $Y_k= V_\ell P_k$ of $A Z_\ell V_{\ell}^H $ with respect to $\range{V_\ell}$, where $Y_k \in \mathbb{C}^{n \times k}$ and
$P_k  = [p_1, \cdots, p_k] \in \mathbb{C}^{\ell \times k}$.
According to Definition \ref{def:hr}, the harmonic Ritz vector $y_j= V_\ell p_j$ then satisfies
\begin{eqnarray} \label{eq:hr:fgmresdr:1}
Z_\ell^H A^H \, (A Z_\ell~ p_j - \theta_j V_\ell~ p_j)  &=& 0. 
\end{eqnarray}
Using the generalized Arnoldi relation (\ref{eq:def:flex_arnoldi}) we finally obtain the relation 
\begin{equation} \label{eq:hr:fgmresdr:2}
{{\bar{H}_\ell}}^H~{{\bar{H}_\ell}}~y_j = \theta {{\bar{H}_\ell}}^H V_{\ell+1}^H V_\ell~y_j.
\end{equation}
Since 
\[
\bar{H}_\ell = \begin{bmatrix}
H_\ell \\
h_{\ell+1,\ell} e_\ell^T\\ 
\end{bmatrix}, \quad H_\ell \in \mathbb{C}^{\ell \times \ell}
\]
where $H_\ell \in \mathbb{C}^{\ell \times \ell}$ is supposed to be nonsingular, the generalized eigenvalue problem is
then equivalent to
\begin{equation} \label{eq:gmres-dr:hr_eig}
(H_\ell + h_{\ell+1,\ell}^2  H_\ell^{-H} e_\ell e_\ell^T)  y_j = \theta_j y_j.
\end{equation}
This corresponds to a standard eigenvalue problem of dimension $\ell$ only, where $\ell$ is supposed to be much smaller 
than the problem dimension $n$. In consequence, the approximate spectral information based on Harmonic Ritz pair is
quite inexpensive to compute.

\paragraph{\label{sec:nonh:aug:app:gmresdr} GMRES augmented with approximate invariant subspace} The augmentation \linebreak space $\cal{W}$ based on approximate invariant information corresponding to $\range{Y_k}$ is then used. The key point detailed next is to understand how to incorporate this information in a minimum residual norm subspace method such as GMRES. To do so, we recall a useful relation satisfied by the harmonic Ritz vectors $P_k \in \mathbb{C}^{\ell \times k}$ shown in \cite[Lemma 3.1]{ggpv:10}
\begin{eqnarray} 
AZ_{\ell} P_k  =   V_{\ell+1}
\begin{bmatrix}
\begin{bmatrix}
       P_k   \\
       0_{1 \times k}
   \end{bmatrix}, c -\bar{H}_\ell y^\star
\end{bmatrix}
\left[
\begin{array}{c}
{\rm diag}(\theta_1,\ldots, \theta_k)\\
\alpha_{1 \times k}
\end{array}
\right], \label{lemma:hr:fgmresdr:eq:expanded} \\
AZ_{\ell} P_k  =   [V_{\ell} P_k,  r_0]
\left[
\begin{array}{c}
{\rm diag}(\theta_1,\ldots, \theta_k)\\
\alpha_{1 \times k}
\end{array}
\right],\label{lemma:hr:fgmresdr:eq}
\end{eqnarray}
where  $r_0 = V_{\ell+1} (c -\bar{H}_\ell y^\star)$ and $\alpha_{1 \times k}=[\alpha_1, \ldots, \alpha_k] \in \mathbb{C}^{1 \times k}$.
Next, the QR factorization of the $(\ell+1) \times (k+1)$ matrix appearing on the right-hand side of relation 
(\ref{lemma:hr:fgmresdr:eq:expanded}) is performed as
\begin{eqnarray}
\label{eq:fgmresdr:qr}
\begin{bmatrix}\begin{bmatrix}\gmresdralg{P_k}\\{0}_{1\times k}\end{bmatrix}, \gmresdralg{c} -\gmresdralg{\bar{H}_\ell} y^\star  \end{bmatrix} =  QR
\end{eqnarray} 
where $Q \in \Complex^{(\ell+1) \times (k+1)}$ has orthonormal columns and $R \in \Complex^{(k+1) \times (k+1)}$ is upper triangular, respectively. Then it can be shown that the relations
\begin{eqnarray} 
A~Z_k &=& V_{k+1} \bar{H}_k, \label{eq:fgmresdr:z}\\
V_{k+1}^H V_{k+1} &=& I_{k+1},\label{eq:fgmresdr:v}\\
\range{[Y_k, r_0]} &=& \range{V_{k+1}}, \label{eq:fgmresdr:range}
\end{eqnarray}
hold with {\it{new}} matrices $Z_k, V_{k} \in\mathbb{C}^{n \times k}$ and $\bar{H}_k\in\mathbb{C}^{(k+1) \times k}$ defined as
\begin{eqnarray} 
Z_k   &=& Z_{\ell}~Q_{\ell\times k},   \label{eq:fgmresdr:z:def}\\
V_{k+1} &=& V_{\ell+1}~Q, \label{eq:fgmresdr:v:def}\\
\bar{H}_k &=& Q^H~\bar{H}_\ell~Q_{\ell\times k}, \label{eq:fgmresdr:h:def}
\end{eqnarray}
where $V_{\ell+1}$, $Z_{\ell}$ and $\bar{H}_\ell$ refer to matrices obtained at the end of the 
previous cycle; see \cite[Proposition 2]{ggpv:10}. With the augmentation subspace ${\cal{W}} = \range{Y_k}$, $m$ Arnoldi steps with variable preconditioners and starting vector $v_{k+1}$ are then carried out while maintaining orthogonality to $V_{k}$  leading to
\begin{equation*}
A~ [z_{k+1}, \cdots, z_{m+k} ] = [v_{k+1}, \cdots, v_{m+k+1} ]~ \bar{H}_{m} \quad \mbox{and}
\quad V_{m+k+1}^H~V_{m+k+1} = I_{m+k+1}.
\end{equation*}
We note that $\bar{H}_{m} \in \mathbb{C}^{(m+1) \times m}$ is upper Hessenberg. At the end of the new cycle this gives the generalized Arnoldi relation
\begin{eqnarray*} 
A ~[Z_{k}, z_{k+1}, \cdots, z_{m+k}] &=& [V_{m+k+1} ]~ 
\begin{bmatrix}
\begin{bmatrix} \bar{H}_k\\{0}_{m\times k}\end{bmatrix} \,  
\begin{bmatrix} B_{k \times m}\\ \bar{H}_{m}\end{bmatrix}
\end{bmatrix} \label{sec:flexible:fgmresdr:bf}
\end{eqnarray*}
i.e.
\[
A Z_{m+k} = V_{m+k+1} \bar{H}_{m+k},
\]
where $V_{m+k+1}\in\mathbb{C}^{n \times (m+k+1)}$, $\bar{H}_{m+k}\in\mathbb{C}^{(m+k+1) \times (m+k)}$ and $B_{k \times m} \in \mathbb{C}^{k \times m}$ results from the orthogonalization of $[Az_{k+1}, \cdots, Az_{m+k+1}]$ against $V_{k+1}$. We note that $\bar{H}_{m+k}$ is no more upper Hessenberg due to the leading dense $(k+1) \times k$ submatrix $\bar{H}_k$. 
It is important to notice that the augmentation space varies at each restart since it is built from the search space available at the end of each previous cycle.
The resulting algorithm can be viewed as an adaptive augmented Krylov subspace method.
We refer the reader to \cite[Sections 2 and 3]{ggpv:10} for the complete derivation of the method and additional comments on its computational cost.

\paragraph{\label{sec:nonh:aug:app:remarks} Remarks and applications}

When the preconditioner is fixed, the previous algorithm proposed by Morgan \cite{morg:02} is known as GMRES with deflated restarting (GMRES-DR). Although the term ``deflated" is used, we note that this algorithm does correspond to a GMRES method with an adaptive augmented basis without any explicit deflated matrix. The success of GMRES-DR  has been demonstrated on many academic examples~\cite{morg:95} and concrete applications such as in lattice QCD \cite{damw:04,frnz:12}, reservoir modeling \cite{akkl:09,kwcs:07} or electromagnetism \cite{ggpv:10}. We refer the reader to~\cite{morg:02,rofi:07} for further comments on the algorithm and computational details. We note that GMRES with deflated restarting is equivalent to other augmented GMRES methods such as GMRES with eigenvectors \cite{morg:95} and implicitly restarted GMRES \cite{morg:00}. Most often the approximate invariant subspace is chosen as the Harmonic Ritz pair corresponding to the smallest in modulus Harmonic Ritz values. Depending on the problem we note however that other specific part of the spectrum of the preconditioned operator can be targeted; see, e.g., \cite[Section 4.2]{ggpv:10} for an application related to a wave propagation problem.

\subsection{\label{sec:nonh:def} Deflated Krylov subspace methods}
We next briefly describe minimal residual Krylov subspace methods based on deflation. We refer the reader to \cite{ggln:11,ggln:12,gutk:12} for a recent excellent overview of deflated  Krylov subspace methods in the Hermitian and non-Hermitian cases, where extensive bibliographical references and historical comments can be found. The general idea of deflation is to split the approximation space into two complementary subspaces such that the projected linear system, referred to as the deflated linear system, will be easier to solve iteratively than the original linear system (\ref{eq:orig}). The fact that these subspaces can be chosen in different ways explains the huge literature on deflated Krylov subspace methods. The Krylov subspace method is then confined in one of this subspace, by projecting the initial residual into this space and by replacing $A$ by its restriction to this space. If the projection operator is chosen properly the deflated linear system will be easier to solve iteratively than the original linear system (\ref{eq:orig}). This property will be notably shown for Hermitian positive 
definite systems in Section \ref{sec:hpd:def} and can be extended to non-Hermitian situations with additional assumptions on $A$ (see, e.g., \cite[Section 2]{erna:08}).  We first present a possible strategy based on orthogonal projection and then briefly discuss an extension based on oblique projection proposed in \cite{gutk:12}.

\subsubsection{\label{sec:nonh:def:orth} Deflation based on orthogonal projection}

We still denote by ${\cal{W}}$ a subspace of $\mathbb{C}^{n}$ of dimension $k$, where $k$ is assumed to be much smaller than the problem dimension $n$. We later denote by $W \in\mathbb{C}^{n \times k}$ a matrix whose columns form a basis of ${\cal{W}}$ so that $W^H A^H A W $ is HPD (hence invertible). 
To simplify further developments, we introduce the matrices $\Qu, \Pu, \Pd \in \mathbb{C}^{n \times n}$ defined respectively as
\begin{eqnarray}
\Qu    &=& A W (W^H A^H A W)^{-1} W^H A^H, \\
\Pu   &=& I_n - \Qu,\\
\Pd &=& I_n - W (W^H A^H A W)^{-1} W^H A^H A.
\end{eqnarray}
We can easily show that $\Pu$ and $\Pd$ are orthogonal projectors such that $\Pu$ projects onto $(A{\cal{W}})^{\perp}$ along $(A{\cal{W}})$, whereas $\Pd$ projects onto ${\cal{W}}^{\perp}$ along ${\cal{W}}$. Furthermore we note that $\Pu$ is Hermitian and that $A \Pd = \Pu A$. The decomposition based on orthogonal projection reads as
\[
\mathbb{C}^{n} = {\cal{W}} \oplus {\cal{W}}^{\perp}.
\]
Hence, the solution $x^\star$ of the original system (\ref{eq:orig}) can be written as
\[
x^\star= (I-\Pd) x^\star + \Pd x^\star = W (W^H A^H A W)^{-1} W^H A^H b + \Pd x^\star.
\]
With this decomposition the original system (\ref{eq:orig}) then simply becomes
 \begin{equation}
\Pu A x = \Pu b. \label{eq:defl}
\end{equation}
Although the deflated matrix $\Pu A$ is singular, the deflated linear system (\ref{eq:defl}) is consistent so that it can 
be solved by an appropriate Krylov subspace method. Here we focus on the application of minimum residual Krylov subspace method based on GMRES to solve the deflated linear system (\ref{eq:defl}). Hence, the search space of the Krylov subspace method applied to (\ref{eq:defl}) can be written as
\[
\hat{S}_m = {\cal{K}}_m(\Pu A, \Pu r_0),
\]
while the current approximation $\hat{x}_m$ and the current residual $\hat{r}_m = \Pu (b - A\hat{x}_m)$ at the end of the cycle satisfies the relations
\begin{eqnarray*}
\hat{x}_m             &\in&   \hat{x}_0 + \hat{S}_m, \\
\Pu (b - A\hat{x}_m)  &\perp& \Pu A~{\cal{K}}_m(\Pu A, \Pu r_0).
\end{eqnarray*}

Since $\Pu A W = 0_{n \times k}$, $\Pu A$ is singular. Hence it is of paramount importance to analyze the possibilities of a breakdown when solving the deflated linear system (\ref{eq:defl}). In our context, when GMRES is used to solve the deflated linear system, this feature has been notably analyzed in \cite[Section 3]{gutk:12} based on theoretical results obtained by Brown and Walker \cite{brwa:97}. We refer the reader to \cite[Corollary 3]{gutk:12} for conditions that characterize the possibility of breakdowns. It is worthwhile to note that a breakdown cannot occur if the condition 
\[
Ker(\Pu A) \cap \range{\Pu A} = \{0 \}
\]
holds; see \cite[Theorem 4.1]{ggln:11}. This condition is notably satisfied if ${\cal{W}}$ is chosen as an exact $A$-invariant subspace, i.e., when $A {\cal{W}} = {\cal{W}}$ since $Ker(\Pu A) = {\cal{W}}$ and $Im(\Pu A) = {\cal{W}}^{\perp}$ due to the nonsingularity of $A$. Once the solution of the deflated linear system is obtained, we deduce the approximation $x_m$ of the original system as
\[
x_m = W (W^H A^H A W)^{-1} W^H A^H b + \Pd \hat{x}_m,
\]
and by construction we note that 
\[
b - A x_m = \Pu (b - A\hat{x}_m),
\]
i.e.,
\[
r_m = \hat{r}_m.
\]

We refer to \cite{erna:08} for applications of deflated Krylov subspace methods with orthogonal projection to linear systems with non-Hermitian matrices. As an illustration, a typical choice of subspaces is to choose the columns of $W$ as right eigenvectors of $A$ corresponding to eigenvalues of small absolute value. 

\subsubsection{\label{sec:nonh:def:obli} Deflation based on oblique projection}

We briefly mention a strategy based on oblique projection that is considered as more appropriate for the solution of non-Hermitian linear systems since the eigenspaces of $A$ are in general not mutually orthogonal \cite{gutk:12}. As in Section \ref{sec:nonh:def:orth}, the search space $S_{\ell}$ will be decomposed into a direct sum of two subspaces. More precisely, the following decompositions into nonorthogonal complements are used
\[
\mathbb{C}^{n} = A{\cal{W}} \oplus \tilde{\cal{W}}^{\perp} = A{\tilde{\cal{W}}} \oplus {\cal{W}}^{\perp},
\]
where ${\cal{W}}$ and $\tilde{\cal{W}}$ represent two subspaces of $\mathbb{C}^{n}$ of dimension $k$ respectively. As before, we denote by $W \in\mathbb{C}^{n \times k}$ ($\tilde{W} \in \mathbb{C}^{n \times k}$) a matrix whose columns form a basis of ${\cal{W}}$ ($\tilde{\cal{W}}$ respectively). We assume that both matrices are chosen such that $\tilde{W}^H A W$ is nonsingular. The key idea is then to introduce the matrices $\Qd, \Pt \in \mathbb{C}^{n \times n}$ defined as
\begin{eqnarray}
\Qd   &=& W (\tilde{W}^H A W)^{-1} \tilde{W}^H , \\
\Pt   &=& I_n - W (\tilde{W}^H A W)^{-1} \tilde{W}^H.
\end{eqnarray} 
It is easy to show that $\Qd$ and $\Pt=I_n-\Qd$ are projection operators; $\Qd$ is the oblique projection onto $(A{\cal{W}})$ along $\tilde{\cal{W}}^{\perp}$, while $\Pt$ is the oblique projection onto $\tilde{\cal{W}}^{\perp} $ along $(A{\cal{W}})$. Given these oblique projection operators, the deflated linear system is now defined as
\[
\Pt A \Pt x = \Pt b
\]
with $\hat{r}_0 = \Pt r_0$. The use of a Krylov space solver is then now restricted to  $\tilde{\cal{W}}^{\perp} $. Hence, 
it can be shown that the deflated Krylov subspace method based on GMRES yields iterates ${x}_m$ at the end of the cycle of the form
\begin{eqnarray*}
{x}_m        &\in & x_0 + {\cal{K}}_m(\Pt A \Pt, \Pt {r}_0) + \cal{W}.
\end{eqnarray*}
This also implies the following relation for the residual \cite{gutk:12}
\begin{eqnarray*}
b  - A{x}_m  &\in & {r}_0 + A {\cal{K}}_m(\Pt A \Pt, \Pt {r}_0) + A{\cal{W}}.
\end{eqnarray*}
We refer the reader to \cite[Sections 5 and 6]{gutk:12} for the mathematical aspects of deflated Krylov subspace methods based on oblique projection and to \cite[Section 11]{gutk:12} for an overview of partly related methods that only differ in the choice of the 
projection operators. A typical choice  is to choose the columns of $W$ as right eigenvectors of $A$ and the columns of $\tilde{W}$ as the corresponding left eigenvectors. We refer to \cite{erna:08} for an application of deflated Krylov subspace methods with oblique projection in the general non-Hermitian case. 

\subsubsection{\label{sec:nonh:def:prec} Deflation by preconditioning}

Finally, we note that deflation based on spectral approximate information can be used to construct nonsingular preconditioners that move small in modulus eigenvalues away from zero. Both Kharchenko and Yeremin \cite{khye:95} and Erhel et al. \cite{erbp:96} have proposed GMRES algorithms with augmented basis and a nonsingular right preconditioner that move the small eigenvalues to a (multiple) large eigenvalue.  Baglama et al.~\cite{bcgr:98} have proposed a left preconditioned GMRES method with similar effect. In \cite{khye:95} the main idea is to translate a group of small eigenvalues of $A$ via low-rank projections of the form
\[
\tilde{A} = A (I_n + u_1 w_1^{H}) \cdots (I_n + u_k w_k^{H}),
\]
where $u_j$ and $w_j$ are the right and left eigenvectors associated with the eigenvalues to be translated respectively. The restarted Krylov subspace method is now applied to the coefficient matrix $\tilde{A}$ leading to an adaptive update of the preconditioner (performed at the end of each cycle). We note that $A$ can correspond to an already preconditioned operator, in such a case this strategy leads to a two-level preconditioning strategy that is found to be effective on real-life applications provided that the spectral information is computed accurately \cite{cadg:03}. We also mention the extension of this two-level preconditioning strategy to the case of sequences of linear systems (see, e.g., \cite{gigm:07} where additional theoretical results and numerical experiments can be found).

\subsection{\label{sec:nonh:augdef} Augmented and deflated Krylov subspace methods}

In the previous sections, we have described how either augmentation or deflation can be incorporated into Krylov subspace methods of minimum residual norm type. We note that it is possible to combine simultaneously deflation and augmentation in a single Krylov subspace method. In such a setting, the search space of the Krylov subspace method is then decomposed as
\[
S_\ell = {\cal{W}} + {\cal{K}}_m(\hat{A}, \hat{r}_0) 
\]
where ${\cal{W}}$ is the augmentation space of dimension $k$, $\hat{A}$ refers to the deflated operator and $\hat{r}_0$ to the deflated residual. As an illustration, we  review the GCRO (Generalized Conjugate Residual with Orthogonalization) method due to de Sturler \cite{stur:96}. 

\subsubsection{\label{sec:nonh:augdef:equi} Equivalence between deflated and augmented methods}

In this section, we describe a general setting that helps us to understand the link between 
deflated and augmented minimal residual norm Krylov subspace methods. It has been first presented in 
\cite{gutk:12} and we generalize this setting to the case of flexible methods. As discussed in Section \ref{sec:nonh:aug}, the search space in augmented methods is of the form
\[
S_\ell = {\cal{W}} \oplus {\cal{K}}_m(\hat{A}, \hat{r}_0)
\]
where ${\cal{W}}$ is an augmentation subspace of dimension $k$. The approximation $x_m$ at the end of a given cycle can be written as
\[
x_m = x_0 + Z_m y_m + W w_m
\]
where $y_m \in \mathbb{C}^{m}$ and $w_m \in  \mathbb{C}^{k}$. In the augmented Krylov subspace methods that we have considered, the residual $r_m$ satisfies a Petrov-Galerkin condition, i.e., $r_m \perp A S$ which leads to the two orthogonality conditions 
\[
r_m \perp A {\cal{W}} \quad \mbox{and}  \quad r_m \perp A {\cal{K}}_m(\hat{A}, \hat{r}_0).
\]
The first orthogonality condition $r_m \perp A {\cal{W}}$ leads to the relation
\[
 (W^H A^H A W) w_m =  W^H A^H (r_0 - AZ_m y_m).
\]
To simplify notation we introduce the matrix $\Qq \in\mathbb{C}^{n \times n}$ such that
\[
\Qq   = W (W^H A^H A W)^{-1} W^H.
\]
We then deduce the following relations for the current approximation $x_m$
\begin{equation}
\label{eq:sec:nonh:equi:x}
x_m = (I_n- \Qq A^H A) (x_0 + Z_m y_m) + \Qq A^H b,
\end{equation}
and for the current residual $r_m$
\begin{equation}
\label{eq:sec:nonh:equi:r}
r_m = (I_n- A \Qq A^H) (r_0 - AZ_m y_m).
\end{equation}
We then introduce the two matrices
\begin{eqnarray*}
\Pc &=& I_n - \Qq A^H A,\\
\Ps &=& I_n - A \Qq A^H
\end{eqnarray*}
where $\Pc  \in\mathbb{C}^{n \times n}$ and $\Ps  \in\mathbb{C}^{n \times n}$. It is easy to show that both $\Pc$ and $\Ps$ are orthogonal projectors and that $A \Pc = \Ps A$. If we define $ \tilde{x}_m \in \mathbb{C}^{n}$ as $\tilde{x}_m = x_0+Z_m y_m$
then relations (\ref{eq:sec:nonh:equi:x}) and (\ref{eq:sec:nonh:equi:r}) become
\begin{eqnarray}
x_m &=& \Pc \tilde{x}_m + \Qq A^H b, \\ 
r_m &=& \Ps (b - A \tilde{x}_m).
\end{eqnarray}
Finally the second orthogonality condition $r_m \perp A {\cal{K}}_m(\hat{A}, \hat{r}_0)$ can then be stated as
\[
r_m = \Ps (b - A \tilde{x}_m) \perp A {\cal{K}}_m(\Ps A, \Ps{r}_0).
\]
We summarize these developments in the following proposition (see \cite[Theorem 2.2]{gutk:12}). 

\begin{Prop}
The following two sets of conditions
\begin{eqnarray*}
x_m &\in& x_0 + {\cal{W}} + {\cal{K}}_m(\hat{A}, \hat{r}_0), \\
r_m &=& b - A x_m \perp ( A {\cal{W}} + A {\cal{K}}_m(\hat{A}, \hat{r}_0) ),
\end{eqnarray*}
and 
\begin{eqnarray*}
\tilde{x}_m &\in& x_0 + {\cal{K}}_m(\hat{A}, \hat{r}_0), \\
\tilde{r}_m &= &\Ps (b - A \tilde{x}_m) \perp  A {\cal{K}}_m(\hat{A}, \hat{r}_0)
\end{eqnarray*}
are equivalent in the sense that
\begin{eqnarray}
x_m &=& \Pc \tilde{x}_m + \Qq A^H b \quad \mbox{and} \quad r_m = \tilde{r}_m.
\end{eqnarray}
\end{Prop}

The first set of conditions corresponds to the standard augmentation approach described in Section \ref{sec:nonh:aug}. In this class of methods the augmentation space $\cal{W}$ is explicitly included in the search space $S$ of the minimum residual Krylov subspace method and $\hat{A} =A$, $\hat{r}_0 ={r}_0$. The second set of conditions corresponds to the standard deflation approach described in Section \ref{sec:nonh:def}. Indeed the 
iteration $\tilde{x}_m$ is first obtained such that the residual  $\Ps (b - A \tilde{x}_m)$ satisfies the Petrov-Galerkin orthogonality
condition. Then a correction is added such that $r_m = \tilde{r}_m$. Both approaches are found to be equivalent. They only differ in the way the augmentation subspace is treated (explicitly or implicitly).

\subsubsection{\label{sec:nonh:augdef:gcro} Methods based on augmentation and deflation}

Methods based on both augmentation and deflation have been introduced recently; see, e.g., \cite{bajm:05,stur:96,stur:99,vavu:94}. We focus here on the Generalized Conjugate Residual with inner Orthogonalization (GCRO) \cite{stur:96}, which combines augmentation and deflation judiciously as detailed next.

GCRO belongs to the family of inner-outer methods \cite[Ch. 12]{axel:94} where the outer iteration is based on the Generalized Conjugate Residual method (GCR), a minimum residual norm Krylov subspace method proposed by Eisenstat, Elman and Schultz \cite{eies:83} while the inner part is based on GMRES respectively. Following the theoretical framework introduced in \cite{eies:00}, GCR maintains a correction subspace spanned by $\range{Z_k}$ and an approximation subspace spanned by $\range{V_k}$, where $Z_k, V_{k} 
\in\mathbb{C}^{n \times k}$ satisfy the relations
 \begin{eqnarray*} 
A~ Z_k &=& V_{k},\\
V_{k}^H ~ V_{k} &=& I_k.
\end{eqnarray*}
 The optimal solution of the minimization problem $ \min \norm{b - A x} $ over the subspace $x_0 + \range{Z_k}$ is then found as $x_k = x_0 + Z_k ~V_{k}^H ~r_0$. Consequently $r_k = b - A x_k$ satisfies
\begin{equation*}
r_k = r_0  - V_{k}~ V_{k}^H r_0 = (I_n - V_{k}~ V_{k}^H) r_0, \quad r_k \perp \range{V_k}.
\end{equation*} 
In \cite{stur:96} de Sturler suggested that the inner iteration takes place in a subspace orthogonal to the outer Krylov subspace. In this inner iteration the following projected linear system is considered 
$$
 (I_n - V_{k}~ V_{k}^H) A z = (I_n - V_{k}~ V_{k}^H) r_k = r_k .
$$
 The inner iteration is then based on a deflated linear system with $(I_n - V_{k}~ V_{k}^H)$ as orthogonal projection.  If a minimum residual norm subspace method is used in the inner iteration to solve this projected linear system approximately, the residual over both the inner and outer subspaces are minimized. Hence, augmentation is applied in the outer iteration and deflation in the inner part of the method. Numerical experiments (see, e.g.,  \cite{stur:96} and \cite[Chapter 1]{fokk:96}) indicate that the resulting method may perform better than other inner-outer methods (without orthogonalization) in some cases. 
 
 We mention that the augmentation subspace can be based on spectral approximate invariant subspace information. This leads to the GCRO with deflated restarting method (GCRO-DR) \cite{psmj:06} that uses Harmonic Ritz information to define the augmentation subspace as in Section \ref{sec:nonh:aug}. This method has been further extended to accommodate variable preconditioning leading to the FGCRO-DR method \cite{cglv:11}. We also refer the 
 reader to \cite{cglv:11} for additional comments on the computational cost of FGCRO-DR and a detailed comparison with the flexible variant of GMRES-DR. When a fixed right preconditioner is used, GMRES-DR and GCRO-DR are equivalent. When variable preconditioning is considered, it is however worthwhile to note that FGMRES-DR and FGCRO-DR are only equivalent if a certain collinearity condition given in \cite[Theorem 3.6]{cglv:11} is satisfied. 
 
  In \cite{stur:99} de Sturler proposed to define an augmentation subspace based on information other than approximate spectral invariant subspace. At the end of each cycle, the strategy (named GCRO with optimal truncation (GCROT)) decides which part of the current global search subspace to keep to define the new augmentation subspace such that the smallest inner residual norm is obtained. This truncation is done by examining angles between subspaces and requires specification of six different parameters that affect the truncation. We refer to \cite{stur:99} for a complete derivation of the method and numerical experiments (see also \cite[Section 4.5]{eies:00}). Finally we note that the extension of GCROT to the case of variable preconditioning has been proposed in \cite{hizi:10} with application to aerodynamics.

\section{\label{sec:hpd} Hermitian positive definite matrices}

Similarly to unsymmetric problems both augmentation and deflation can be considered to speed-up the convergence of the conjugate gradient method,
possibly in combination with preconditioning.
However, contrarily to the previous methods based on Arnoldi basis construction, the conjugate gradient method relies on a short term 
recurrence and restarting mechanisms do not need to be implemented to control the memory consumption.
Consequently the space used for augmentation or for deflation should be fully defined before starting the iteration for a given right-hand side.

\subsection{\label{sec:hpd:aug} Augmented conjugate gradient methods}
As discussed in Section~\ref{sec:nonh:aug}, the search space in augmented methods 
$S_\ell = W \oplus K_m(A,r_0)$ is a $\ell$ dimensional space (with $\ell = m +k$) where ${\cal W}$ is an augmentation subspace of dimension $k$ spanned by $k$ linearly
independent vectors $W = [ w_1, \cdots w_k]$.
In order to build a basis of this space a deflated Lanczos algorithm can be used, that consists in applying a standard
Lanczos method starting from an unit  vector $v_1$ using the matrix
$$
B = A - A W (W^H A W)^{-1} W^H A.
$$
Notice that $W^H A W$ is full rank since $A$ is HPD. If $v_1$ is orthogonal to $W$, deflated Lanczos  builds a sequence of orthonormal vectors $V_m= [ v_1, \cdots, v_m]$ ($V^H_m A V_m = I_m$) that spans a space orthogonal to $W$, i.e., $W^H V_m = 0_{k \times m}$.
For the Krylov subpace part of $S_\ell$, $x_0$ is chosen so that $r_0 = b- Ax_0 \perp W$ and $v_1 = r_0/\norm{r_0}$.
That can be guaranteed by defining  $x_0$ from any $x_{-1}$ as $ x_0 = x_{-1} + W (W^H A W)^{-1} W^H r_{-1}$.
The augmented CG algorithm seeks for a solution $x_\ell = x_0 +  W \mu_{\ell} + V_m y_\ell \in x_0 +  W + K_m(A,r_0)$ with the Ritz-Galerkin condition
$r_\ell = b - A x_\ell \perp ( W + K_m(A,r_0))$.
Using the above described space and orthogonality condition, it is shown in~\cite{syeg:00}, that the following properties (that are very similar and inherited from the classical CG) still hold.
\begin{Prop}
The iterate $x_j$, the residual $r_j$ and the descent directions $p_j$ satisfy the following relations and properties
\begin{itemize}
  \item $r_j$ is collinear to $v_{j+1}$, that is, the residual vectors are orthogonal to each other,
  \item the short term recurrences are satisfied:
    \begin{eqnarray*}
       x_j & = & x_{j-1} + \alpha_{j-1} p_{j -1} \\
       r_j & = & r_{j-1} - \alpha_{j-1} A p_{j -1} \\
       p_j & = & r_{j} + \beta_{j-1} p_{j-1} - W \mu_j 
    \end{eqnarray*}
    where $\alpha_{j -1}$ and $\beta_{j -1}$ have the same expression as in classical CG and $$\mu_j = (W^H A W)^{-1}W^H A r_j,$$
   \item the vectors  $p_j$ are A-orthogonal to each other as well as A-orthogonal to all the $w_j$'s.
\end{itemize}
\end{Prop}

Using theoretical results from~\cite{ergu:00}, the following properties related to convergence rate and
optimization property of the iterate are shown in~\cite{syeg:00}. 

\begin{Prop}
The approximate solution $x_\ell$ is such that
\begin{itemize}
 \item  the convergence history exhibits an upper bound expression on the convergence rate similar to classical CG
 \begin{equation}\label{eq:augCGCvg}
\displaystyle \norm{x_\ell - x^\star}_A \leq 
2\cdot \left ( \frac{ \sqrt{\kappa(\Hp^HA\Hp)} - 1}{\sqrt{\kappa(\Hp^HA\Hp)} + 1} \right)^\ell
    \norm{x_0 - x^\star}_A, 
\end{equation}
where $\kappa( \cdot)$ denotes the condition number of the matrix and $\Hp$ is the A-orthogonal projection on $W^{\perp_A}$.
This projection is defined by $\Hp = I_n - W (W^H A W)^{-1} W^H A$.
 \item  similarly to classical CG, the iterate complies with a minimum A-norm error on the search space
   $\displaystyle x_\ell = \argmin_{x \in x_0 +  W + K_m(A,r_0) } \norm{x - x^\star}_A$.
\end{itemize}
\end{Prop}

\paragraph{Augmenting using an invariant subspace}
Let $(\lambda_1=\lambda_{\min}, \cdots, \lambda_s=\lambda_{\max})$ denote the $s$ distinct eigenvalues of $A$ ordered by increasing magnitude (i.e., values as they are real positive). The invariant subspace spanned by the $k$ extreme (either largest or smallest) eigenvalues can be
used in place of ${\cal W}$ to build the augmented space. 
Equation~\eqref{eq:augCGCvg} shows  that $\displaystyle \kappa(A) = \frac{\lambda_{\max}}{\lambda_{\min}}$ (that would appear in this bound for classical
CG) is replaced either by $\displaystyle \frac{\lambda_{\max}}{\lambda_{k+1}}$ if the left most part of the spectrum is used or
by  $\displaystyle \frac{\lambda_{s-k}}{\lambda_{\min}}$ if the right most part is used. 
Consequently if $\lambda_{\min} \ll \lambda_{k+1}$ ($\lambda_{s-k} \ll \lambda_{\max}$) the convergence of augmented CG should be significantly faster
than the convergence of CG on the original system.

\subsection{\label{sec:hpd:def} Deflated Krylov subspace methods}
We next briefly describe  CG variants based on deflation. 
As mentioned in Section~\ref{sec:nonh:def}, the general idea of deflation is to split the approximation space into two complementary subspaces.
Similarly to the notation in the previous sections, we denote by ${\cal{W}}$ a subspace of $\mathbb{C}^{n}$ of dimension $k$, where $k$ is
 assumed to be much smaller than the problem dimension $n$. 
We later denote by $W \in\mathbb{C}^{n \times k}$ a matrix whose columns form a basis of ${\cal{W}}$.
Because $A$ is Hermitian positive definite, $W^H A W$ is also HPD and hence invertible.
We can then define the following projector
\begin{equation}\label{eq:projectHPD}
 \Pq = I - W(W^H A W)^{-1}W^HA
\end{equation}
that is an oblique projector along $W$ ($\Pq$ is equal to $\Hp$).
As in the non-Hermitian case, we decompose the solution $x^\star =(I-\Pq)x^\star+\Pq x^\star$ and compute each component separately.
In particular, $(I-\Pq)x^\star = W(W^HAW)^{-1}W^HAx^\star=W(W^HAW)^{-1}W^H b$ essentially reduces to the solution of a 
small $k \times k$ system.
For the calculation of the second component $\Pq x^\star$, it can be observed that $A\Pq=\Pq^HA$ so that $A\Pq x^\star = \Pq^HA x^\star = \Pq^H b$.
Even though the matrix $\Pq^HA$ is Hermitian semi-definite positive of rank
$n-k$ (its nullspace is $W$), CG can still be used because the deflated linear system  $\Pq^HA x^\star = \Pq^H b$ is consistent~\cite{kaas:88}.
Furthermore,  because the null space never enters the iteration, the corresponding zero eigenvalues do not influence the convergence~\cite{kaas:88}
and we can define the effective condition number of the positive semidefinite matrix $\Pq^HA$, denoted $ \kappa_{eff}(\Pq^HA)$, as
the ratio of its largest to smallest strictly positive eigenvalues.

Once the linear system $\Pq^H A \tilde{x} = \Pq^H b$ is solved, one just needs to apply $\Pq$ to this solution to
compute the second component of the solution.
This technique still requires the solution of a linear system of size $n$ using the CG method, but is expected to be
more effective if $\kappa_{eff}(\Pq^HA) \ll \kappa(A)$. We refer the reader to \cite{frvu:01} for a discussion on the choice of
$W$.

\paragraph{Deflating using an invariant subspace}
If $W$ defines an invariant subspace of $A$ associated with extreme eigenvalues, the situation becomes much clearer.

Let assume that $W$ defines an invariant subspace associated with the smallest eigenvalues 
$(\lambda_1, ..., \lambda_k)$ of $A$. We have $\Pq^HAW=0_{m \times k}$ so that $\Pq^HA$ has $k$ zero eigenvalues. Because $A$ is HPD,  $Z = W^{\perp}$, the  orthogonal complement of $W$ (i.e., $W^H Z = 0$ so that $\Pq^H Z = Z$) defines
an invariant subspace associated with the eigenvalues $\lambda_{k+1}, ..., \lambda_n = \lambda_{\max}$.
Therefore, we have $A Z = Z B $ for some nonsingular $B$.
Consequently we have $\Pq^HAZ=\Pq^HZB=ZB$ so that $Z$ is an invariant subspace of $\Pq^H A$ associated
with the same eigenvalues $\lambda_{k+1}, ..., \lambda_{\max}$.
This shows that
$$
\kappa_{eff}(\Pq^HA) = \frac{\lambda_{\max}}{\lambda_{k+1}},
$$
that indicates that deflating using an invariant subspace cancels the corresponding eigenvalues, leaving the rest of the spectrum unchanged.
If $ \lambda_{k+1} \gg \lambda_1 = \lambda_{\min}$  the convergence of CG is significantly speeded-up.

\subsection{\label{sec:hpd:defprec} Deflation via preconditioning}
Using spectral information, it is possible to design preconditioners that enable to exhibit a condition
number for the preconditioned matrix similar to $\kappa_{eff}(\Pq^HA)$.

Let $W=[w_1, \dots w_k]\in \Bbb{C}^{n \times k}$ be the  normalized eigenvectors of $A$ associated with $\{\lambda_i\}_{i=1, \dots, k}$ the set of smallest eigenvalues. Let  $\nu$ be a real positive value. We can then define the preconditioner
$$
 \MDeflat = I_n + W ( \nu (W^H A W)^{-1} - I_k ) W^H.
$$
This preconditioner is such that $\MDeflat A W= \nu W$ and $\MDeflat Aw=Aw$ if 
$W^Hw=0$ (in particular any eigenvectors of $A$ not in $W$), which shows that $\MDeflat$ moves the eigenvalues $\{\lambda_i\}_{i=1, \dots, k}$ 
to $\nu$ and leaves the rest of the spectrum unchanged. 
If $\nu = \lambda_{k+1}$, the condition number of the preconditioned matrix is the same as the one of the
deflated matrix in the previous section.

Furthermore we can define additive coarse space correction preconditioners inspired from domain decomposition techniques.
They lead to preconditioned matrices with similar condition number as well.
We then define
$$
  \MShift = I_n + \nu W (W^H A W)^{-1} W^H.
$$
This preconditioner is such that $\MShift Aw_i=(\nu +\lambda_i)w_i$ and $\MShift Aw=Aw$ if 
$W^Hw=0$.
That is, the eigenvalues $\{\lambda_i\}_{i=1, \dots, k}$ 
are shifted to $\nu+\lambda_i$, while the rest of the spectrum is unchanged. 
If it exists $\nu$ so that $\lambda_{k+1} \leq \lambda_{\min} + \nu \leq  \lambda_k + \nu \leq \lambda_{\max}$,
the preconditioned matrix would have again the same condition number as the one of the deflated system $\kappa_{eff}$.

We refer the reader to~\cite{gigr:06} for an analysis of the condition number of this class of preconditioners when
approximated spectral information is used.
We also refer to~\cite{tnve:09} and the references therein for the exposure of various preconditioning techniques 
that can be defined using various combinations of these building box components.


\section{\label{sec:nonh:mrhs} Linear systems with multiple right-hand sides given in sequence}

Although our primary focus is the solution of a single linear system with preconditioned Krylov subspace methods, it is however possible to include deflation and augmentation in a broader setting. Indeed in many applications in computational science and engineering, linear systems with multiple right-hand sides have to be solved. More precisely we are interested in solving a sequence of linear systems defined as $A^{l} x^{l} = b^{l}$ where both the non-Hermitian matrix $A^{l} \in \mathbb{C}^{n \times n}$ and the right-hand side $b^{l}\in \mathbb{C}^{n}$ may change from one system to the next, and the linear systems may typically not be available simultaneously. If we consider a sequence of identical or {\it{slowly}} changing matrices, Krylov subspace methods based on augmentation and deflation are appropriate since subspace recycling is then possible. The key idea is to extract relevant information (e.g. approximate invariant subspace but not only) while solving a given system, and then to use this information to further accelerate the convergence of the Krylov subspace method for the subsequent linear systems. At this point, augmented or deflated Krylov subspace methods of Sections \ref{sec:nonh:aug}, \ref{sec:nonh:def} and \ref{sec:nonh:augdef} can then be used. We refer the reader to \cite[Chapter 3]{park:05} for a detailed analysis of subspace recycling in the non-Hermitian case and to \cite{psmj:06} where the GCRO method augmented with approximate spectral information is shown to be efficient on applications related to fatigue and fracture of engineering components, electronic structure calculation and quantum chromodynamics; see also \cite{kist:06} for an application in optical tomography. Recent applications are related to model reduction \cite{befe:11} (see also \cite{asgc:12} for recycling methods based on BiCG).

For HPD matrices, if a sequence of linear systems with the same matrix but different right-hand sides has to be solved
different alternatives can be considered to define the space to augment the search space from one solve to the next.
In~\cite{ergu:00}, an approach based on harmonic Ritz values is described that might be implemented using only the
first $m \geq k$ steps of augmented CG iteration.
Still to reduce the memory footprint of the eigenvector calculation, in~\cite{stor:10} a thick-restart Lanczos is embedded in the CG 
iterations to extract accurate spectral information.

\section{\label{sec:conc} Conclusions and prospectives}

We have briefly reviewed the main features and mathematical properties of augmented and deflated Krylov subspace methods 
for the solution of certain linear systems of equations where the coefficient matrix was either non-Hermitian or Hermitian positive
definite. These increasingly popular procedures combined with preconditioning have been shown effective on a wide range of applications in computational science and engineering as mentioned in this paper. We are certainly aware that this brief overview is far from being complete. Results related to two-sided Krylov subspace methods in the non-Hermitian case or the treatment of the Hermitian indefinite case are indeed missing; see, e.g., \cite{asgc:12,ggln:11,ggln:12,wasp:07} for additional comments and references. Similarly, the solution of linear systems with multiple right-hand sides given at once has not been covered. For such a class of problems, augmented and deflated block Krylov subspace methods have been studied (see, e.g., \cite{morg:05,yusj:02}) and their efficiency has been proved on realistic applications. Finally we would like to mention that algebraic connections between deflation, multigrid and domain decomposition have been made in recent papers \cite{klrh:12,navu:08,tnve:09}.

Concerning implementation aspects, some augmentation and deflation procedures are already present 
in the main software projects such as either PETSc\footnote{http://www.mcs.anl.gov/petsc/} or Trilinos\footnote{http://trilinos.sandia.gov/} for the solution of large-scale, complex multi-physics engineering and 
scientific problems. More precisely, in its scalable linear equation solvers (KSP) component, PETSc includes an algorithm described in \cite{erbp:96}, while the Belos package in Trilinos notably proposes an augmented and deflated approach based on GCRO-DR \cite{psmj:06}. Most likely there will be a growing effort to incorporate augmented and deflated Krylov subspace methods in such libraries in a near future. Finally designing variants or new Krylov subspace methods for the next generation of massively parallel computing platforms is currently a topic of active research in the numerical linear algebra community;
see \cite{gamv:13,piplinedCG,hoem:10} for algorithms, comments and references. Thus in a near future it is highly probable that variants of augmented and deflated Krylov subspace methods will be proposed as well.


%
%
\topsep=0.0ex
\parsep=0.0ex
\parskip=0.0ex
\itemsep=0.0ex

%
%
\nocite{*}

\bibliographystyle{plain}

\end{document}